\documentclass[12pt]{article}
\usepackage{pstricks,pst-node}
\usepackage{pst-plot}

\usepackage{amsmath,amsthm}
\usepackage{amssymb}
\usepackage{color}
\usepackage{xspace}
\usepackage[colorlinks=true,
linkcolor=green,
filecolor=brown,
citecolor=green]{hyperref}
\def\red{\textcolor{red} }
\def\blue{\textcolor{blue} }

\def\d{\ensuremath{\mathcal D}\xspace}
\def\k{\ensuremath{\mathcal K}\xspace}

\begin{document}
\newtheorem{theorem}{Theorem}
\newtheorem{defn}[theorem]{Definition}
\newtheorem{lemma}[theorem]{Lemma}
\newtheorem{prop}[theorem]{Proposition}
\newtheorem{cor}[theorem]{Corollary}
\begin{center}
{\Large
A bijection for Delannoy paths          \\[2mm]
David Callan  \\
}
\end{center}
 
\begin{abstract}
We exhibit a bijection between central Delannoy $n$-paths, that is, lattice paths from the origin to $(n,n)$ with steps $E=(1,0), \,N=(0,1),\,D=(1,1)$ and the lattice paths from 
the origin to $(n+1,n)$ where the only restriction on the steps is that they have finite nonnegative slope. 
\end{abstract}


\section{Introduction} \label{intro}

A \emph{Delannoy path} \cite{delannoy2005} is a lattice path consisting of steps $E=(1,0), \,N=(0,1),\,D=(1,1)$. For $m,n \ge 0$, let $\d_{m,n}$ denote the set of Delannoy paths from $(0,0)$ to $(m,n)$ and set $\d_n=\d_{n,n}$, the \emph{central Delannoy paths} \cite{sulanke03}. Clearly, a central Delannoy path has the same number of $E$ and $N$ steps. A path in $\d_n$ with $k$ East ($E$) steps is determined by a linear arrangement of $k$ $E$s, $k$ $N$s, and $n-k$ $D$s. Thus there are $\binom{n+k}{k,k,n-k}=\binom{n}{k}\binom{n+k}{k}$ paths in $\d_n$ with $k$ $E$ steps, sequence \htmladdnormallink{A063007}{http://oeis.org/A063007} in OEIS \cite{oeis}.

A \emph{Kimberling path}, after sequence \htmladdnormallink{A049600}{http://oeis.org/A049600}, is a lattice path from the origin where the only restriction on the steps is that they have finite nonnegative slope (vertical steps not allowed). Let $\k_{i,j}$ denote the Kimberling paths that terminate at $(i,j)$. The 3 paths in $\k_{2,1}$ are shown in Figure 1. 
\begin{center}
\begin{pspicture}(-5,-2)(5,1.5)
\psline[linecolor=blue](-6,0)(-5,0)(-4,1)
\psline[linecolor=blue](-1,0)(0,1)(1,1)
\psline[linecolor=blue](4,0)(6,1)
 
\psdots(-5,0)(-6,0)(-4,1)(-1,0)(0,1)(1,1)(4,0)(6,1)
\psset{linecolor=lightgray}
\psdots(-4,0)(0,0)(1,0)(5,0)(6,0)
\psdots(-6,1)(-5,1)(-1,1)(4,1)(5,1)

\rput(-6,-.3){\textrm{{\footnotesize $O$}}}
\rput(-3.5,1.1){\textrm{{\footnotesize (2,1)}}}
\rput(-1,-.3){\textrm{{\footnotesize $O$}}}
\rput(1.5,1.1){\textrm{{\footnotesize (2,1)}}}
\rput(4,-.3){\textrm{{\footnotesize $O$}}}
\rput(6.5,1.1){\textrm{{\footnotesize (2,1)}}}

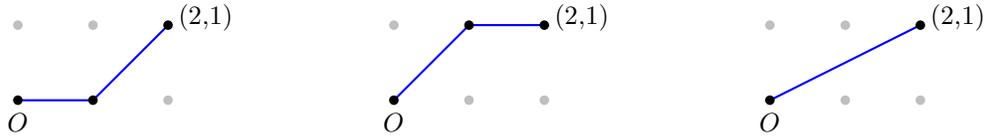
\rput(0,-1.2){\textrm{{\small The paths in $\k_{2,1}$: the middle path has vertex set $\{(0,0),(1,1),(2,1)\}$}}} 

\rput(0,-2){\textrm{{\small Figure 1}}}  

\end{pspicture}
\end{center}  
Suppose a Kimberling path in $\k_{i,j}$ has $k$  interior vertices. Their $x$-coordinates form a $k$-set (distinct elements) from $\{1,\dots,i-1\}$ and, independently, their $y$-coordinates form a $j$-multiset (repetition allowed) from $\{0,1,\dots,j\}$. Hence, the number of paths in $\k_{i,j}$ with $k$ interior vertices is $\binom{i-1}{k}\binom{j+k}{k}$. In particular, $\k_{n+1,n}$ contains $\binom{n}{k}\binom{n+k}{k}$ paths with $k$ interior vertices.

In the next section, we give a bijection between the equinumerous sets of the previous two paragraphs.

\section{A bijection $\mathbf{\phi:\ \d_\emph{n} \to \k_{\emph{n}\,+1,\emph{n}}}$} \label{bij}

Given a central Delannoy path, such as $P=NEEDNNNEDDEEN \in \d_n$ with $n=9$  $k=5$ East steps, label each East and North step with the $y$-coordinate of its terminal point as in Figure 2a. 

\begin{pspicture}(0,-2.5)(15,6)
\psset{
  xunit = .6, yunit = .6,
  algebraic
}
  \psaxes[dx =1,dy = 1,labels = none,
    ticklinestyle = dotted,tickwidth = .1pt,
    xticksize = 0 8,yticksize = 0 8
  ]{}(0,0)(0,0)(8,8)[$x$,0][,0]
  \psaxes[labels = none,
    dx = 1,Dx = 1,dy = 1,Dy = 1
  ](0,0)(0,0)(8,8)
  \psaxes[
    dx =1,
    dy = 1,
    labels = none,
    ticklinestyle = dotted,
    tickwidth = 0.5pt,
    xticksize = 0 8,
    yticksize = 0 9
  ]{}(15,0)(15,0)(24,8)[$x$,0][,0]
 \psaxes[labels = none,
    dx = 1,Dx = 1,dy = 1,Dy = 1
  ](15,0)(15,0)(24,8)

\psline[linestyle=solid,linecolor=gray](0,0)(8,8)
\psline[linestyle=solid,linecolor=gray](15,0)(24,8)

\psline[linecolor = blue,linewidth=1.4pt](0,0)(0,1)(1,1)(2,1)(3,2)(3,3)(3,4)(3,5)(4,5)(5,6)(6,7)(7,7)(8,7)(8,8)

\psline[linecolor = blue,linewidth=1.4pt](15,0)(16,1)(18,1)(19,5)(20,7)(23,7)(24,8)
\psset{linecolor=black,linewidth=1.2pt}
\psdots(15,0)(16,1)(18,1)(19,5)(20,7)(23,7)(24,8)
\psset{linecolor=black,linewidth=.8pt}
\psdots(0,0)(0,1)(1,1)(2,1)(3,2)(3,3)(3,4)(3,5)(4,5)(5,6)(6,7)(7,7)(8,7)(8,8)

\rput(11,5){$\longrightarrow$}
\rput(0,8.5){$y$}
\rput(15,8.5){$y$}

\red{\rput(.3,.5){\textrm{1}}
\rput(2.7,2.5){\textrm{3}}
\rput(2.7,3.5){\textrm{4}}
\rput(2.7,4.5){\textrm{5}}
\rput(7.7,7.6){\textrm{8}} }

\blue{
\rput(.3,1.4){\textrm{1}}
\rput(1.3,1.4){\textrm{1}}
\rput(3.3,5.4){\textrm{5}}
\rput(6.3,6.6){\textrm{7}}
\rput(7.3,6.6){\textrm{7}}}

\rput(4,-1.8){central Delannoy path $P$}
\rput(20,-1.8){Kimberling path $\phi(P)$}

\rput(4,-3.5){Figure 2a}
\rput(20,-3.5){Figure 2b}

\rput(.8,-.6){\footnotesize{1}}
\rput(1.8,-.6){\footnotesize{2}}
\rput(2.8,-.6){\footnotesize{3}}
\rput(3.8,-.6){\footnotesize{4}}
\rput(4.8,-.6){\footnotesize{5}}
\rput(5.8,-.6){\footnotesize{6}}
\rput(6.8,-.6){\footnotesize{7}}
\rput(7.8,-.6){\footnotesize{8}}

\rput(-.7,1){\footnotesize{1}}
\rput(-.7,2){\footnotesize{2}}
\rput(-.7,3){\footnotesize{3}}
\rput(-.7,4){\footnotesize{4}}
\rput(-.7,5){\footnotesize{5}}
\rput(-.7,6){\footnotesize{6}}
\rput(-.7,7){\footnotesize{7}}
\rput(-.7,8){\footnotesize{8}}

\rput(15.8,-.6){\footnotesize{1}}
\rput(16.8,-.6){\footnotesize{2}}
\rput(17.8,-.6){\footnotesize{3}}
\rput(18.8,-.6){\footnotesize{4}}
\rput(19.8,-.6){\footnotesize{5}}
\rput(20.8,-.6){\footnotesize{6}}
\rput(21.8,-.6){\footnotesize{7}}
\rput(22.8,-.6){\footnotesize{8}}
\rput(23.8,-.6){\footnotesize{9}}

\rput(14.3,1){\footnotesize{1}}
\rput(14.3,2){\footnotesize{2}}
\rput(14.3,3){\footnotesize{3}}
\rput(14.3,4){\footnotesize{4}}
\rput(14.3,5){\footnotesize{5}}
\rput(14.3,6){\footnotesize{6}}
\rput(14.3,7){\footnotesize{7}}
\rput(14.3,8){\footnotesize{8}}

\end{pspicture}

The image path $\phi(P)$ in $\k_{n+1,n}$ will have $k$ interior vertices and the $i$-th interior vertex is given by $(x_i,y_i)$ where $x_i$ and $y_i$ are the labels on the $i$-th North and East steps respectively, $1\le i \le k$.  In the example, the interior vertices, left to right, are $(\red{1}, \blue{1}), (\red{3}, \blue{1}), (\red{4}, \blue{5}), (\red{5}, \blue{7}), (\red{8}, \blue{7})$ as in Figure 2b. The line joining the initial and terminal vertices is shown for each path in Figure 2 in anticipation of the next section. 

To invert $\phi$, let $A$ be the set of $x$-coordinates and $B$ the multiset of $y$-coordinates of the interior vertices of a Kimberling path in $\k_{n+1,n}$. Set $C=\{1,\dots,n\}\,\backslash\, A$. In the example illustrated in Figure 2, we have $A=\{1,3,4,5,8\},$ $B=\{1,1,5,7,7\},\ C=\{2,6,7\}$. Then $A,B,C$ give the heights above the $x$-axis of the $N,E,D$ steps respectively in the order specified as follows. First arrange the (disjoint) union of $A$ and $B$ in weakly increasing order with the $A$ number first in case of equality: 
$1^A\,1^B\,1^B\,3^A\,4^A\,5^A\,5^B\,7^B\,7^B\,8^A$. Then insert the $C$ numbers into this list so that each is as far left as possible while maintaining the weakly increasing property:  $1^A\,1^B\,1^B\,2^C\,3^A\,4^A\,5^A\,5^B\,6^C\,7^C\,7^B\,7^B\,8^A$. This yields the original central Delannoy path by changing $A$ to $N$, $B$ to $E$, and $C$ to $D$.

\section{An application of  $\mathbf{\phi}$} \label{props}
A \emph{subdiagonal} lattice path is one that lies weakly below the line joining its initial and terminal vertices. Subdiagonal central Delannoy paths are known as Schr\"{o}der paths, counted by the large Schr\"{o}der numbers, \htmladdnormallink{A006318}{http://oeis.org/A006318}.
\begin{prop}
The bijection $\phi:\d_n \to \k_{n+1,n}$ preserves the subdiagonal property.
\end{prop}
\begin{proof}
We will establish a stronger result: the $i$th East step in $P$ lies (weakly) above the diagonal $y=x$ if and only if the $i$-th interior vertex $(x_i,y_i)$ of $\phi(P)$ lies above its diagonal line $y=\frac{n}{n+1}x$ joining (0,0) and $(n+1,n)$. 
To see  this, for $P\in \d_n$ with $k$ East steps, let $(X_i,Y_i),\ 1\le i \le k,$ denote the $(x,y)$ coordinates of 
the east end of the $i$-th East step. Thus the $i$-th East step in $P$ is above $y=x$ iff $Y_i/X_i \ge 1$, 
and $ (x_i,y_i)$ is above $y=\frac{n}{n+1}x$ iff $y_i/x_i >\frac{n}{n+1}$. (Note that, by integer divisibility considerations, $y_i/x_i$ never equals $\frac{n}{n+1}$.)

So suppose $u$ $D$s precede the $i$-th $N$ in $P$ and $v$ $D$s precede the $i$-th $E$. Thus $x_i=i+u$ and $X_i=i+v$. Also, $Y_i=y_i$ by the definitions.
\vspace*{-5mm}
\begin{itemize}
\item If $u=v$, then $X_i=i+u =x_i$ and the condition $Y_i/X_i \ge 1$ implies $  y_i/x_i =Y_i/X_i \ge 1>\frac{n}{n+1}$.  Conversely, 
$y_i/x_i>\frac{n}{n+1} $ implies $ y_i/x_i\ge 1 $ because $\frac{n}{n+1} <y_i/x_i<1$ with $x_i,y_i$ positive integers would imply $y_i \le x_i-1$ and hence $x_i>n+1$.
\item If $u<v$, then $Y_i\ge i+v$ because $v>u$ implies there is a $D$ after the $i$th  $N$ and before the $i$th $E$. So $Y_i \ge X_i$  and   $Y_i / X_i\ge 1$. Also, $y_i/x_i \ge \frac{i+v}{i+u} >1 > \frac{n}{n+1}$.

\item If $u>v$,  then $Y_i < i+v$ because $v<u$ implies the $i$th  $E$ occurs before the $i$th $N$ and so fewer than $i$ $N$s precede the $i$th $E$. Hence $Y_i/X_i <1$. Also, $y_i=Y_i\le i+v-1$ and $x_i=i+u>i+v\ge i+v+1$ and so $y_i/x_i \le \frac{i+v-1}{i+v+1}<\frac{n}{n+1}$ since $i+v<2n+1$.
\end{itemize} 
\end{proof}
\textbf{Corollary.} \ Subdiagonal Kimberling paths in $\k_{n+1,n}$ are counted by the $n$-th large Schr\"{o}der number.

Department of Statistics, University of Wisconsin-Madison

\end{document}